\documentclass[11pt]{article}

\usepackage[numbers,sort&compress]{natbib}
\usepackage{hyperref}
\usepackage{amsmath}
\usepackage{amssymb}
\usepackage{amsfonts}
\usepackage{amsthm}

\def\init{\setcounter{equation}{0}}

\newtheorem{theorem}{Theorem}[section]

\newtheorem{proposition}[theorem]{Proposition}

\newtheorem{lemma}[theorem]{Lemma}

\newtheorem{definition}[theorem]{Definition}

\def\eps{\epsilon}

\def\zz{\mathbb{Z}_N}
\def\ee{\mathbb{E}}

\def\pp{\mathbb{P}}

\def\cc{\mathbb{C}}
\def\nn{\mathbb{N}}
\def\rr{\mathbb{R}}

\begin{document}

\title{Arithmetic structures in random sets}
\author{Mariah Hamel and Izabella {\L}aba}

\maketitle

\begin{abstract}

We extend two well-known results in additive number theory, S\'ark\"ozy's theorem on square differences in dense sets and a theorem of Green
on long arithmetic progressions in sumsets, to subsets of random sets of asymptotic
density 0.  Our proofs rely on a restriction-type Fourier analytic argument of Green 
and Green-Tao.

\end{abstract}


\section{Introduction}

\label{Introduction}

\init

The purpose of this paper is to extend several basic results
in additive number theory, known for sets of positive density in $\zz$,
to the setting of random sets of asymptotic density 0.  
This line of work originated in the paper of Kohayakawa-{\L}uczak-R\"odl
\cite{KLR}, who proved a random-set analogue of Roth's theorem
on 3-term arithmetic progressions. Roth's theorem \cite{roth} asserts that for any fixed 
$\delta>0$ there is a large integer $N_0$ such that if $N>N_0$
and if $A$ is a subset of $\{1,\dots,N\}$ with $|A|\geq\delta N$,
then $A$ contains a non-trivial 3-term arithmetic progression
$a,a+r,a+2r$ with $r\neq 0$. The article \cite{KLR} raises the following 
question: are there any sets $W$, sparse in $\{1,\dots,N\}$,
with the property that any set $A$ containing a positive proportion
of the elements of $W$ must contain a 3-term arithmetic progression?
The authors proceed to answer it in the affirmative for random sets:

\begin{theorem}\label{klr-theorem} \cite{KLR}
Suppose that $W$ is a random subset of $\zz$ such that the events
$x\in W$, where $x$ ranges over $\zz$, are independent and have
probability $p=p(N)\in(CN^{-1/2},1]$. Fix $\alpha>0$. Then the
statement
\begin{quote}
every set $A\subset W$ with $|A|\geq \alpha |W|$ contains a 3-term
arithmetic progression
\end{quote}
is true with probability $1-o_\alpha(1)$ as $N\to\infty$.

\end{theorem}

The current interest in questions of this type is motivated by the work of 
Green \cite{green-roth} and Green-Tao \cite{GT05}, \cite{GT06}
on arithmetic progressions in the primes, where
the ``pseudorandomness" of the almost-primes plays a key role.
For example, Tao-Vu \cite[Section 10.2]{TV} give
an alternative (and simpler) proof of Theorem \ref{klr-theorem} under the
stronger assumption that $p\geq CN^{-\theta}$ with $\theta$ small
enough.  While the argument in
\cite{KLR} is combinatorial and uses Szemer\'edi's regularity lemma,
the proof in \cite{TV} is Fourier-analytic and relies in particular
on a restriction-type estimate from \cite{green-roth}, \cite{GT06}.

It is natural to ask which other results from additive number theory can
be extended to the random set setting.  While the methods of \cite{KLR}
do not seem to extend to other questions, the decomposition 
technique in \cite{GT06} turns out to be more robust.  We are able to
use it to prove random set analogues of two well-known results:
S\'ark\"ozy's theorem on square differences, and a theorem of Green
on long arithmetic progressions in sumsets.  

We now give the precise statement of our results.  Throughout the paper,
$W$ is a random subset of $\zz$, with each $x\in\zz$
belonging to $W$ independently with probability $p\in(0,1]$.  We will
assume that $p\geq N^{-\theta}$, where $\theta$ is a sufficiently
small positive number.  In particular, we allow $p$ to go to 0 as
$N\to\infty$. We also fix $\delta>0$ and let $A\subset W$,
$|A|=\delta |W|$.

S\'ark\"ozy's theorem (proved also independently by Furstenberg)
states that for any fixed positive
number $\delta$ there is a large integer $N_0$ such that if $N>N_0$
and if $A$ is a subset of $\{1,\dots,N\}$ with $|A|\geq\delta N$,
then $A$ contains two distinct elements $x,y$ such that
$x-y$ is a perfect square.  The best known quantitative bound,
due to Pintz, Steiger and Szemer\'{e}di \cite{PSS}, is that
one may take $N_0=(\log N)^{-c\log\log\log\log
N}$.  In the converse direction, Ruzsa
\cite{R84} constructed a set of size $N^{1-0.267}$ which contains no
square difference.

We are able to prove the following.

\begin{theorem}\label{sarkozyextension} Suppose that $W$ is a random
subset of $\zz$ such that the events $x\in W$, where $x$ ranges over
$\zz$, are independent and have probability
$p=p(N)\in(cN^{-\theta},1]$ where $0<\theta<1/110$. Let $\alpha>0$.
Then the statement
    \begin{quote} for every set $A\subset W$ with $|A|\geq \alpha
    W$, there are $x,y\in A$ such that $x-y$ is a non-zero perfect square
    \end{quote}
is true with probability $o_{\alpha}(1)$ as $N\rightarrow\infty$.
\end{theorem}

We also have an analogous result for higher power differences, see Section
\ref{powerdifferences}.

If $A,B$ are two sets of
integers, we will write $A+B=\{a+b:\ a\in A,b\in B\}$.
Let $W$ be a random set as described above, but with $\theta\in(1/2,1]$.
One can show using a probabilistic argument that
it holds with probability $1-o(1)$ that the sumset $A+A$ of every subset $A\subset W$
with $|A|>\alpha |W|$ has density at least $\alpha^2$ in $\zz$\footnote{
We are grateful to Mihalis Kolountzakis for pointing this out to us and communicating
a short proof.}.  If $\theta$ is close enough to 1, then we can prove the following stronger
result using Fourier-analytic methods.

\begin{proposition}\label{sumsets-size}
Suppose that $W$ is a random subset of $\zz$ such that the events
$x\in W$, where $x$ ranges over $\zz$, are independent and have
probability $p=p(N)\in(CN^{-\theta},1]$, where $0<\theta<1/140$. Then
for every $\beta<\alpha$, the statement
\begin{quote}
for every set $A\subset W$ with $|A|\geq \alpha |W|$,
we have $|A+A|\geq \beta N$
\end{quote}
is true with probability $1-o_{\alpha,\beta}(1)$ as $N\to\infty$.

\end{proposition}

It is easy to see that one can have $|A+A|\approx \alpha N$ in the setting of the
proposition: let $A_x=W\cap(P+x)$, where $P$ is an arithmetic progression in $\zz$
of step about $\alpha^{-1}$ and length about $\alpha N$.  An averaging argument shows
that $|A_x|\gg\alpha|W|$ for some $x$, while $|A_x+A_x| \leq 2|P|\approx\alpha N$.

Our second main result concerns the existence
of long arithmetic progressions in sumsets.  Bourgain
\cite{B90} proved that if $A,B$ are sumsets of $\{1,\dots,N\}$ with
$|A|>\alpha N$, $|B|>\beta N$, then $A+B$ contains a $k$-term
arithmetic progression with
\begin{equation}\label{ap-e1}
k>\exp(c(\alpha\beta\log N)^{1/3}-\log\log N).
\end{equation}
The point here is that a sumset has much more arithmetic structure, and therefore
contains much longer arithmetic progressions, than would be normally expected
in a set of a similar size (based on Szemer\'edi's theorem, for example).
Bourgain's bound was improved by Green \cite{green-sumsets} to
\begin{equation}\label{ap-e2}
k>\exp(c(\alpha\beta\log N)^{1/2}-\log\log N),
\end{equation}
which is the best known result in this direction so far. An alternative proof of
essentially the same bound was given more recently by Sanders \cite{S06}.
On the
other hand, Ruzsa \cite{R91} gave a construction showing that the
exponent $1/2$ in (\ref{ap-e2}) cannot be improved beyond $2/3$.
Note that if $A=B$, the estimate (\ref{ap-e2}) gives a non-trivial
result only when $\alpha>(\log N)^{-1/2}$, and in particular sets
with density $N^{-\epsilon}$ are not allowed.

The case of sparse sets was considered more recently by
Croot-Ruzsa-Schoen \cite{CRS}.  The authors proved that if
$A,B\subset\zz$ obey $|A||B|\geq(6N)^{2-\frac{2}{k-1}}$, then $A+B$
contains a $k$-term arithmetic progression.  They also gave a
construction of sets $A\subset\zz$ with $|A|\geq N^{1-\theta}$,
where $\theta$ is small enough depending on $\epsilon
>0$, such that $A+A$
does not contain an arithmetic progression longer than
$\exp(c\theta^{-\frac{2}{3} -\epsilon})$.

Our result is the following.

\begin{theorem}\label{sumsets-theorem}
Suppose that $W$ is a random subset of $\zz$ such that the events
$x\in W$, where $x$ ranges over $\zz$, are independent and have
probability $p=p(N)\in(CN^{-\theta},1]$, where $0<\theta<1/140$. Assume that
$\alpha$ and $k$ obey
\begin{equation}\label{e-alpha}
\alpha\geq\frac{C_1\log\log N}{\sqrt{\log N}},
\end{equation}
\begin{equation}\label{e-k}
k\leq \exp\left( \frac{\alpha^2\log\log N}{C_2\log\frac{1}{\alpha}(\log\log\log N+\log\frac{1}{\alpha})}
\right),
\end{equation}
where $C_1,C_2$ are sufficiently large constants.  Then the statement
\begin{quote}
for every set $A\subset W$ with $|A|\geq \alpha |W|$, the sumset
$A+A$ contains a k-term arithmetic progression
\end{quote}
is true with probability $1-o_{k,\alpha}(1)$ as $N\to\infty$.

\end{theorem}

A non-quantitative version of the result, namely that the displayed statement in
the theorem is true with
probability $1-o(1)$ as $N\to\infty$ if $\alpha$ and $k$ are fixed, can be obtained
by applying Szemer\'edi's theorem to the positive density set $A+A$.
Our point, as in \cite{B90} or \cite{green-sumsets},
is that the arithmetic progressions indicated by Theorem \ref{sumsets-theorem} are
much longer than those in Szemer\'edi's theorem, and that they can be found using a much
easier argument.  For comparison, the current best bounds in Szemer\'edi's theorem
\cite{gowers} imply that a set of relative density $\alpha$ in $\zz$ should contain
$k$-term arithmetic progressions with
$$
k\leq\log\log\left(\frac{\log\log N}{\log\frac{1}{\alpha}}\right),
$$
which is much weaker than (\ref{e-k}).

The bounds on $\theta$ in Theorems \ref{sarkozyextension} and \ref{sumsets-theorem}
are due to our choices of exponents in the proofs and are probably not optimal.
The natural threshold would be $1/2$, as in \cite{KLR}.
However, it does not seem possible to extend our results to all $\theta>1/2$ using
the same type of arguments as in this paper.

The article is organized as follows.  In the next section we explain the notation and
summarize the known results that will be used repeatedly.
Theorem \ref{sarkozyextension} is proved in Sections \ref{varnavidessection} and
\ref{sarkozy}.  Its analogue for higher power differences, Theorem \ref{powerextension},
is stated and proved in Section \ref{powerdifferences}.
The proof of Theorem \ref{sumsets-theorem} is given in Section \ref{long-section},
with the proofs of the main estimates postponed to Sections \ref{sec-iteration}
and \ref{sec-random}.  The proof of Proposition \ref{sumsets-size}, which
involves a simplified version of the argument in the proof of Theorem
\ref{sumsets-theorem}, concludes the paper.


\section{Preliminaries}
\label{prelim} \init


We first explain the notation.  We use $|A|$ to denote the
cardinality of a set $A\subset\zz$.  The {\em probability} of a set
$A$ is $\pp(A)=N^{-1}|A|$, and the {\em expectation} of a function
$f:\zz\to\cc$ is defined as
$$\ee f=\ee_x f=N^{-1}\sum_{x\in\zz}f(x).$$
We will also sometimes use conditional probability and expectation
$$
\pp(A|X)=\frac{|A\cap X|}{|X|},\ \ee(f|X)=\ee_{x\in X}f(x)
=\frac{1}{|X|}\sum_{x\in X}f(x).$$ Whenever the range of a variable
(in a sum, expectation, etc.) is not indicated, it is assumed to be
all of $\zz$. We will also use the notation
$\|f\|_p=(\sum_x|f(x)|^p)^{1/p}$ and $\|f\|_{L^p(X)}=(\sum_{x\in
X}|f(x)|^p)^{1/p}$.  All constants throughout the paper will be independent
of $N$, $\alpha$, and $k$.

The discrete Fourier transform of $f$ is defined by
$$\widehat{f}(\xi)=\ee_x f(x)e^{-2\pi ix\xi/N}.$$
We have the usual Plancherel identity
$\sum\widehat{f}\bar{\widehat{g}}=N^{-1}\sum f\bar{g}$ and the
inversion formula $f(x)=\sum_{\xi\in\zz}\widehat{f}(\xi)  e^{-2\pi
ix\xi/N}$.

We define the convolution of two functions $f,g:\zz\to\cc$ by the
formula
$$(f*g)(x)=\sum_y f(y)g(x-y)=\sum_{t,s:t+s=x}f(t)g(s).$$
We have the identity $N\widehat{f}\widehat{g}=\widehat{f\ast g}$.

We recall a few basic results about Bohr sets, all of which are
standard in the literature and can be found e.g. in
\cite{GT-inverse}, \cite{TV}, or in \cite{B99} where regular Bohr
sets were first introduced.

\begin{definition}\label{def-bohr}
A {\em Bohr set} is a set of the form $B=b+B(\Lambda,\delta)$, where
$b\in\zz$, $\Lambda\subset\zz$, $\delta\in(0,2)$, and
$$B(\Lambda,\delta)=\{x\in\zz:\ |e^{2\pi ix\xi}-1|\leq \delta\hbox{ for all }
\xi\in \Lambda\}.$$
We will often
refer to $|\Lambda|$ and $\delta$ as the rank and radius of
$B$, respectively.
\end{definition}

\begin{definition}\label{def-regular}
Let $c_0$ be a small positive constant which will remain fixed
throughout the paper.  We will say that a Bohr set $B(\Lambda,\delta)$ is {\em regular} if
$$\pp(B(\Lambda,(1+c_0^2)\delta)\setminus B(\Lambda,(1-c_0^2)\delta))
\leq c_0\pp(B(\Lambda,\delta)).$$
We will also say that $B=b+B(\Lambda,\delta)$ is regular if $B(\Lambda,\delta)$ is regular.
\end{definition}

\begin{lemma}\label{lemma-bohrsize}
If $B=B(\Lambda,\delta)$ is a regular Bohr set, then $\pp(B)\geq
(cc_0^2\delta)^{|\Lambda|}$.
\end{lemma}

\begin{lemma}\label{lemma-regular}
Assume that $c_0$ is small enough.  Then for any $\Lambda
\subset\zz$ with $|\Lambda |\leq \sqrt{c_0}N$ and any $\delta_0>0$ there
is a $\delta\in(\frac{\delta_0}{2},\delta_0)$ such that $B(\Lambda ,\delta)$ is regular.
\end{lemma}

We will need a Fourier-analytic argument which first appeared in
\cite{green-roth} in a slightly different formulation and in
\cite{GT06} as stated, and was adapted in \cite{TV} to a random set
setting. Specifically, \cite{green-roth} and \cite{GT06} introduced
the decomposition $f=f_1+f_2$ defined below, where $f_1$ is the
``structured" bounded part, and $f_2$ is unbounded but random. We
will need several results concerning the properties of $f_1$ and
$f_2$, which we collect in the next two lemmas. The first one is
contained in the proofs of \cite[Proposition 5.1]{GT06} or
\cite[Theorem 10.20]{TV}.

\begin{lemma}\label{lemma-decomposition}
Assume that $f:\zz\rightarrow [0,\infty)$ satisfies $\ee(f)\geq
\delta>0$ and
\begin{equation}\label{e-restriction}
\|\widehat{f}\|_{q}\leq M
\end{equation}
for some $2<q<3$. Assume also that $f\leq \nu$, where
$\nu:\zz\rightarrow [0,\infty)$ obeys the pseudorandom condition
\begin{equation}\label{pseudorandomcondition}
\|\hat{\nu}(\xi)-{\bf 1}_{\xi=0}\|_\infty\leq \eta
\end{equation}
for some $0<\eta\leq 1$. Let
$$f_1(x)=\ee(f(x+y_1-y_2):\ y_1,y_2\in B_0),$$
where
$$B_0=\{x:\ |e^{-2\pi i\xi x/N}-1|\leq\epsilon_0,\ \xi\in\Lambda_0\},
\ \Lambda_0=\{\xi:\ |\widehat{f}(\xi)|\geq \epsilon_0\}$$ for some
$\epsilon_0$ to be fixed later.  Let also $f_2(x)=f(x) -f_1(x)$.
Then

(i) $0\leq f_1\leq 1+(1+\pp(B_0)^{-1})\eta,$

(ii) $\ee f_1=\ee f$,

(iii) $\|\widehat{f_2}(\xi)\|_\infty\leq 3(1+\eta)\epsilon_0,$

(iv) $|\widehat{f_i}(\xi)|\leq |\widehat{f}(\xi)|$ for all
$\xi\in\zz$ and $i=1,2$.  In particular, (\ref{e-restriction}) holds
with $f$ replaced by $f_2$.

\end{lemma}

In order to be able to apply Lemma \ref{lemma-decomposition}, we need
to have the estimate (\ref{e-restriction}) for some $2<q<3$.  To
this end we have the following result, based on the Stein-Tomas
argument as used in \cite{green-roth}, \cite{GT06}, and contained in
the form we need in \cite[Lemma 10.22 and proof of Theorem
10.18]{TV}.

\begin{lemma}\label{lemma-restriction}
Let $f$ and $\nu$ be as in Lemma \ref{lemma-decomposition}, except
that instead of (\ref{e-restriction}) we assume that
$$\|\widehat{f}\|_2\leq C\eta^{-\epsilon/4}$$
for some $\epsilon>0$.  Then (\ref{e-restriction}) holds with
$q=2+\epsilon$.

\end{lemma}

We adapt this argument to the random setting as in \cite[Section
10.2]{TV}. Suppose that $W$ is a random subset of $\zz$ such that
each $x\in\zz$ belongs to $W$ independently with probability
$p\in(0,1)$.  We will assume that $p\geq N^{-\theta}$, where
$0<\theta<1/100$. We also fix $\delta>0$ and let $A\subset W$,
$|A|=\delta |W|$. We let
$$\nu(x)=p^{-1}W(x), \ f(x)=p^{-1}A(x).$$

\begin{lemma}\label{lemma-random}
Let $\nu$ and $f$ be the random variables defined above.  Then

\smallskip
(i) $\|\widehat{\nu(\xi)}-{\bf 1}_{\xi=0}\|_\infty=O(N^{-1/5})$ with
probability $1-o(1)$,

\smallskip
(ii) $\|\widehat{f}\|_2^{2} =N^{-1}\|f\|_2^{2} =O(p^{-1})\leq
N^{\theta}$ with probability $1-o(1)$.
\end{lemma}

Part (i) of the lemma follows from well-known probabilistic
arguments. It can be found e.g. in \cite[Corollary 1.9 and Lemma
4.15]{TV}, or extracted from the proof of Lemma 14 in
\cite{green-sumsets}. Observe in particular that (i) with $\xi=0$
says that $\pp(W)=p(1+O(N^{-1/5}))$ with probability $1-o(1)$.  Part
(ii) follows from this and the Plancherel identity.



\section{A Varnavides-type theorem for square differences}

\label{varnavidessection}

\init


The purpose of this section is to prove the following theorem.

\begin{theorem}\label{varnavidesforsquares}  Let $0<\delta\leq 1$
and $N\geq 1$ be a prime integer.  Let $f:\zz\rightarrow [0,1]$ be a
bounded function such that $$\ee f\geq \delta.$$  Then we have
$$\ee(f(n)f(n+r^2)|n,r\in\zz,\ 1\leq r\leq \lfloor
\sqrt{N/3}\rfloor)\geq c(\delta)-o_{\delta}(1).$$ \end{theorem}

Theorem \ref{varnavidesforsquares} strengthens S\'ark\"ozy's theorem
(stated in the introduction)
in the same way in which a theorem of Varnavides \cite{V59} strengthens
Roth's theorem on 3-term arithmetic progressions.  It guarantees the
existence of ``many" square differences in a set of positive density,
instead of just one.

\begin{proof}  The proof combines S\'ark\"ozy's theorem
with a modification of Varnavides's combinatorial argument \cite{V59}.
We first note that it suffices to prove the result for
characteristic functions.  To see this, let $f$ be as in the
theorem, and define $A:=\{n\in\zz:f(n)\geq \delta/2\}$. Then
$|A|\geq \delta N/2$ and
$f\geq \frac{\delta}{2}$ on $A$. Hence,
assuming the result for characteristic functions, we have
$$\ee(f(n)f(n+r^2))\geq \frac{\delta^2}{4}\ee(A(n)A(n+r^2))\geq
\frac{\delta^2}{4}c(\delta/2).$$

We now turn to the proof of the result for characteristic functions.
Let $A\subset \zz$ such that $|A|\geq \delta N$ and $N$ is
sufficiently large. We will consider arithmetic progressions
$P_{x,r}$, given by
\begin{equation}\label{AP} 1 \leq x< x+r^2<\dots<x+(k-1)r^2\leq
N\end{equation}
where $x,r\in\zz$, $r\leq \sqrt{3N}$, and where $k\in\nn$ is chosen
so that the conclusion of S\'ark\"ozy's theorem holds for subsets
of $\{1,...,k\}$ which
have size at least $\frac{1}{2}\delta k$.

Suppose that
\begin{equation}\label{vv-er}
r^2<\frac{\delta N}{k^2}.
\end{equation}
We say that a progression
$P_{x,r}(N)$ as in (\ref{AP}) is \textit{good} if
\begin{equation}\label{gooddensity}|P_{x,r}(N)\cap A|\geq
\frac{1}{2}\delta k.\end{equation}
Let $G_{r}(N)$ denote the set
of good progressions $P_{x,r}(N)$ for a fixed $r$.  We claim that
\begin{equation}\label{sizeofgoodprogression}
|G_{r}(N)|>\frac{1}{4}\delta N.\end{equation}
Indeed, we have
$$|A\cap(kr^2,N-kr^2)|\geq|A|-2kr^2
\geq \delta N-2kr^2\geq \delta(1-\frac{2}{k})N,
$$
where at the last step we used (\ref{vv-er}).
Each $a\in A\cap(kr^2,N-kr^2)$ is contained in exactly $k$
progressions $P(x,r)$.  Hence
$$
\sum_{x:1\leq x< x+(k-1)r^2\leq N}|A\cap P_{x,r}(N)|
\geq k\delta(1-\frac{2}{k})N>\frac{3}{4}\delta kN\ \ \  (k>8).
$$
On the other hand, the number of progressions $P(x,r)$ for a
fixed $r$ is clearly bounded by $N$, hence we have an upper
bound
$$\sum_{x:1\leq x<x+(k-1)r^2\leq N}|A\cap
P_{x,r}(N)|<N\cdot\frac{1}{2}\delta k+G_{r}(N)k.$$  Combining
these bounds yields (\ref{sizeofgoodprogression}) as claimed.

Let $G(N):=\sum_{r: r^2<\frac{\delta N}{k^2}}G_{r}(N)$.  Then
\begin{equation}\label{vv.100}
G(N)\geq \frac{\sqrt{\delta N}}{k}\frac{\delta
N}{4}=c_1(\delta)N^{3/2},
\end{equation}
since $k$ depends only on $\delta$.

By S\'ark\"ozy's theorem, each good progression
$P_{x,r}(N)$ contains a square difference.  We now count the number
of good progressions which may contain a fixed square difference
pair $x,x+r^2$.
Clearly, $x,x+r^2$ can be contained in at
most $k-1$ progressions with step size $r^2$ and at most
$\frac{1}{2}k(k-1)$ progressions with step size $r^2/t$
for integers $t>1$.  Since $k$ depends only on $\delta$, the
total number of progressions containing $x,x+r^2$ is bounded
by $c_2(\delta)$.
Thus the total number of
square differences in $A$ must be at least
$$\frac{c_1(\delta)}{c_2(\delta)}N^{3/2}=c(\delta)N^{3/2}.$$
Subtracting off the trivial progressions (with $r^2=0$) gives the
desired result.

\end{proof}


\section{Proof of Theorem \ref{sarkozyextension}}

\label{sarkozy}

\init


Let $W,A$ be as in Theorem
\ref{sarkozyextension}.  At least one of the sets $A_1=A\cap[0,N/3)$, $A_2=A\cap[N/3,2N/3)$,
$A_3=A\cap[2N/3,N)$, say $A_1$ (the other two cases are identical), has size at least $|A|/3$.
Define $\nu,f$ as in Lemma
\ref{lemma-random}, but with $A$ replaced by $A_1$.
By Lemma
\ref{lemma-random}, the assumptions of Lemma \ref{lemma-restriction}
with $\eta=N^{-1/5}$ and $\epsilon=1/11$ are satisfied with probability $1-o(1)$,
thus (\ref{e-restriction}) holds with $q=23/11$. We will henceforth condition on
these events.
Let $f=f_1+f_2$ as in
Lemma \ref{lemma-decomposition}, with $\epsilon_0=
\epsilon_0(\alpha)$ small enough to be fixed later.  We would like to ensure
that
\begin{equation}\label{pp.e3}
\|f_1\|_\infty\leq 2.
\end{equation}
By Lemma \ref{lemma-decomposition}, this will follow if
\begin{equation}\label{pp.e2}
N^{-1/5}(1+\pp(B_0))<1.
\end{equation}
By Lemma \ref{lemma-bohrsize}, we can estimate
$\pp(B_0)\gg(c\epsilon_0)^{|\Lambda_0|}$, while by
(\ref{e-restriction}) and Chebyshev's inequality we have $|\Lambda_0|\leq
(M/\epsilon_0)^{23/11}$.  Now a short calculation shows that if
\begin{equation}\label{pp.e1}
\log\frac{1}{\epsilon_0}<c_1\log\log N
\end{equation}
with $c_1$ small enough, which we will assume henceforth,
then (\ref{pp.e2}) and (\ref{pp.e3}) hold.

It suffices to prove that
\begin{equation}\label{transference}
\ee(f(x)f(x+r^2)|x,r\in\zz, 1\leq r\leq
\sqrt{N/3})\geq c(\delta)-o_{\delta}(1).
\end{equation}
Indeed, since $A_1\subset[0,N/3)$, any square difference $a-a'=r^2$ with
$a,a'\in A_1$ and $1\leq r^2\leq{N/3}$ must be an actual square difference in
$\mathbb{Z}$, not just a square difference mod $N$.

We write
$f(x)f(x+r^2)=\sum_{i,j=1}^2 f_i(x)f_j(x+r^2)$, and estimate the expectation of
each term.
Applying Theorem \ref{varnavidesforsquares}
to $f_1$, we get a lower bound on the main term
\begin{equation}\label{sa-main}
\ee(f_1(x)f_1(x+r^2)|x,r\in\zz, 1\leq r\leq
\sqrt{N/3})\geq c_1(\delta)-o_{\delta}(1),
\end{equation}
if $N$ is large enough so that (\ref{pp.e1}) holds.
We now turn to the error estimates.
We write
\begin{equation}\label{sa-e2}
\ee(f_2(x)f_2(x+r^2)|x,r\in\zz, 1\leq r\leq
\sqrt{N/3})
=\sqrt{3N}\,\ee(f_2(x)f_2(x+t)S(t)|x,t\in\zz),
\end{equation}
where
$S(\cdot)$ denotes the characteristic function of the squares less than $N/3$.
From Green \cite{G02} we have
the estimate
$$\|\hat{S}\|_{12}\leq 2^{19/12}N^{-1/2},$$
based on a number theoretic
bound on the number of representations of an integer as the sum of six
squares.  Using also Parseval's identity and H\"{o}lder's
inequality, we have
\begin{align*}\ee(f_2(x)&f_2(x+t)S(t)|x,t\in\zz)\\
&=\sum_{\xi\in\zz}|\hat{f_2}(\xi)|^2|\hat{S}(\xi)|\\
    &\leq
    \big(\sum_{\xi\in\zz}|\hat{S}(\xi)|^{12}\big)^{1/12}
    (\sum_{\xi\in\zz}|\hat{f_2}(\xi)|^{24/11})^{11/12}\\
    &\leq 2^{19/12}N^{-1/2}\|\hat{f_2}\|_{23/11}^{23/12}\| \hat{f_2}\|_{\infty}^{1/12}\\
    &\leq CN^{-1/2} \eps_0^{1/12}.
    \end{align*}
Plugging this into (\ref{sa-e2}), we see that
$$\ee(f_2(x)f_2(x+r^2)|x,r\in\zz, 1\leq r\leq
\sqrt{N/3})\leq c_1(\delta)/4$$ if $\eps_0$ was chosen sufficiently
small depending on $\delta$. The ``mixed" error terms are estimated
similarly. Combining the error estimates with (\ref{sa-main}) yields
(\ref{transference}) as desired.


\section{Power differences}
\label{powerdifferences}  \init


In this section we show that a modification of the proof of Theorem
\ref{sarkozyextension} yields an analogous result for higher power differences.

\begin{theorem}\label{powerextension} Suppose that $W$ is a random
subset of $\zz$ such that the events $x\in W$, where $x$ ranges over
$\zz$, are independent and have probability
$p=p(N)\in(cN^{-\theta},1]$ with $0<\theta<\theta_k$, where $\theta_k$
is small enough depending on $k\in\mathbb{N}$. Let $\alpha>0$.
Then the statement
    \begin{quote} for every set $A\subset W$ with $|A|\geq \alpha
    W$, $\exists$ $x,y\in A$ such that $x-y=n^k$ for some $n\in\mathbb{N}$\end{quote}
is true with probability $o_{k,\alpha}(1)$ as $N\rightarrow\infty$.
\end{theorem}

Since the proof is very similar to that of Theorem \ref{sarkozyextension},
we only sketch the main steps. 
Instead  of Theorem \ref{varnavidesforsquares}, we will need a similar 
result for higher powers, which can be proved by exactly the same argument.

\begin{theorem}\label{varnavidesforpowers}  Let $0<\delta\leq 1$,
and let $N\geq 1$ be a prime integer.  Let $f:\zz\rightarrow [0,1]$ be a
bounded function such that $\ee f\geq \delta.$  Then we have
$$\ee(f(n)f(n+r^k)|n,r\in\zz,\ 1\leq r\leq \lfloor
\sqrt[k]{N/3}\rfloor)\geq c(\delta)-o_{\delta}(1).$$ \end{theorem}

We now follow the argument in Section \ref{sarkozy}.  Define $\nu,f,f_1,f_2$ as in
the proof of Theorem \ref{sarkozyextension}.  Applying Theorem
\ref{varnavidesforpowers} to $f_1$, we see that
$$\ee(f_1(x)f_1(x+r^k)|x,r\in\zz, 1\leq r\leq
\sqrt[k]{N/3})\geq c(\delta)-o_{\delta,\epsilon_0, M}(\eta).$$  
To
estimate the error terms, we invoke the asymptotic formula
for Waring's problem (see e.g.\ \cite{Nathanson}), which implies that
\begin{align*}R_{k,3k}(x)&:=|\{(a_1,...,a_{3k})\in\zz|a_1^k+...+a_{3k}^k\equiv x
            \mod N\}
    \leq cN^2.\end{align*}
By convolution and Parseval identities, this translates to
$$\|\widehat{P_k}\|_{6k}\leq c_1 N^{1/k-1},$$ 
where $P_k$ denotes the characteristic function set of $k$-th powers smaller than $N/3$,
and $c,c_1$ are
constants depending on $k$.  Now we are able to estimate the error
terms as in Section \ref{sarkozy}, for example we have
\begin{align*}\ee(f_2(x)f_2(x+r)P_k(r))&\leq\|\hat{P_k}\|_{6k}
    \|\hat{f_2}\|_{(12k-1)/(6k-1)}^{(12k-2)/(12k)}\| \hat{f_2}\|_{\infty}^{1/6k}\\
    &\leq c_1 C N^{1/k-1}\epsilon_0^{1/6k}.
    \end{align*}
At the last step we used that (\ref{e-restriction}) holds with $q=\frac{12k-1}{6k-1}$ if $\theta_k$ 
is small enough.
The proof is finished as in Section \ref{sarkozy}.


\section{Long arithmetic progressions in sumsets}
\label{long-section} \init


We now turn to Theorem \ref{sumsets-theorem}.  In this section we
prove the theorem, modulo the two main estimates
(\ref{vv.e11}), (\ref{vv.e10}) which will be proved in the next
two sections.

Our proof will combine the arguments of Sanders \cite{S06} with
those of Green-Tao \cite{GT06}. Let $W,A$ be as in Theorem
\ref{sumsets-theorem}, and define $\nu,f$ as in Lemma
\ref{lemma-random}. We will show that, with high probability, there
is a reasonably large Bohr set $B$ on which we have $f*f(x)>0$ for
all but a few values of $x$.  But $f*f$ is supported on $A+A$, hence
all but a small fraction of $B$ is contained in $A+A$.  The proof is
concluded by invoking a pigeonholing argument from \cite{S06}, which
says that the portion of $B$ contained in $A+A$ contains a long
arithmetic progression.

The details are as follows.  Fix $k$ (the length of the
progression), and let $\sigma=(16k)^{-1}$.  We will also assume that $k>k_0$ and
$\alpha<\alpha_0$ , where $k_0\in\nn$ is a sufficiently large absolute constant
and $\alpha_0>0$ is a sufficiently small absolute constant.

By Lemma
\ref{lemma-random}, the assumptions of Lemma \ref{lemma-restriction}
with $\eta=N^{-1/5}$ and $\epsilon=1/9$ are satisfied with probability $1-o(1)$,
thus (\ref{e-restriction}) holds with $q=19/9$.
Let $f=f_1+f_2$ as in
Lemma \ref{lemma-decomposition}, with $\epsilon_0=
\epsilon_0(\alpha,\sigma)$ small enough to be fixed later.
We will assume that (\ref{pp.e1}) holds with $c_1$ sufficiently small;
as in Section \ref{sarkozy}, it follows that $\|f_1\|_\infty\leq 2$.

We
need an extension of a result of Sanders \cite{S06}: there are
regular Bohr sets $B:=b+B(\Gamma,\delta)$ and $B':=b+B(\Gamma,\delta')$
such that
\begin{equation}\label{vv.e11}
\Big|\{x\in B':\ f_1*f_1(x)\geq
\frac{\alpha^2}{2}|B|\}\Big|>(1-\sigma)|B'|,
\end{equation}
 and
\begin{equation}\label{pp.e4}
\delta'\gg\frac{\alpha^2\delta}{|\Gamma|},
\end{equation}
\begin{equation}\label{pp.e5}
\delta\gg\big(\frac{\alpha}{\log(\sigma^{-1})}\big)^{C\log(\alpha^{-1})},
\end{equation}
\begin{equation}\label{pp.e6}
|\Gamma|\ll\alpha^{-2}\log(\sigma^{-1}).
\end{equation}
We establish this in Proposition \ref{f1estimate}. We then verify in
Section \ref{sec-random}, via a restriction-type argument, that if
\begin{equation}\label{pp.e8}
\log\frac{1}{\epsilon_0}\gg \alpha^{-2}\log\frac{1}{\alpha}\log k
(\log\log k+\log\frac{1}{\alpha}),
\end{equation}
with a large enough implicit constant, then
\begin{equation}\label{e-random}
\Big|\{x\in B':\ |f_2*f_i(x)|\geq
\frac{\alpha^2}{10}|B|\}\Big|<\sigma|B'|,\ i=1,2.
\end{equation}
It follows that
\begin{equation}\label{vv.e10}
\Big|\{x\in B':\ f*f(x)\geq
\frac{\alpha^2}{10}|B|\}\Big|>(1-4\sigma)|B'|,
\end{equation}
provided that both (\ref{pp.e1}) and (\ref{pp.e8}) hold.  A somewhat cumbersome
calculation shows that $\epsilon_0$ can be chosen so as to satisfy both (\ref{pp.e1})
and (\ref{pp.e8}), provided that
\begin{equation}\label{pp.e9}
\log k\ll\frac{\alpha^2\log\log N}{\log\frac{1}{\alpha}(\log\log\log N+\log\frac{1}{\alpha})},
\end{equation}
which is equivalent to (\ref{e-k}).

We now invoke Lemma 6.5 in \cite{S06}, which says that if
\begin{equation}\label{pp.e10}
(4\sigma)^{-1}\ll |\Gamma|^{-1}\delta'N^{1/|\Gamma|},
\end{equation}
then the
set on the left side of (\ref{vv.e10}) contains an arithmetic
progression of length $(16\sigma)^{-1}=k$.
Plugging in (\ref{pp.e4})--(\ref{pp.e6}) and solving for $N$, we see that (\ref{pp.e10})
holds if
\begin{equation}\label{pp.e11}
\log N\gg\alpha^{-2}(\log^2 k+\log^2(\frac{1}{\alpha})+\log\frac{1}{\alpha}\log\log k).
\end{equation}
Another cumbersome calculation shows that  if we assume (\ref{pp.e9}), then the additional
condition (\ref{e-alpha}) suffices to guarantee that
(\ref{pp.e11}) holds.
Thus, assuming both (\ref{e-alpha}) and (\ref{e-k}),
the set on the left side of (\ref{vv.e10})
contains a $k$-term arithmetic progression. Since that set is
contained in $A+A$, the conclusion of the theorem follows.

In the next two sections we complete the proof by verifying the
inequalities (\ref{vv.e11}), (\ref{e-random}).


\section{The main term estimate}

\label{sec-iteration}

\init


\begin{proposition}\label{iteration}  Let $B=b+B(\Gamma,\delta)$ be a regular Bohr set.  Let
$f: \zz\rightarrow \rr$ be a function such that
$\textrm{supp}(f)\subset B$, $0\leq f\leq 1$ and
$\ee_{B}f=\alpha>0$. Fix $\sigma\in(0,1]$ and let $d=|\Gamma|$. Then
one of the following must be true:

(i) There is a $\delta'\gg\frac{\alpha^2\delta}{d}$  such that
$B'=b+B(\Gamma,\delta')$ is regular and
\begin{equation}\label{pospropbohr}
\Big|\{x\in B':(f*f)(x)\geq \frac{\alpha^2}{2}|B|\} \Big|\geq
(1-\sigma)|B'|,\end{equation} or

(ii) There is a regular Bohr set
$B''=b''+B(\Gamma\cup\Lambda,\delta'')$ such that
\begin{equation}\label{densityincrement}
\ee(f|B'')\geq\alpha\big(1+2^{-5}\big),
\end{equation}
where $|\Lambda|\ll \alpha^{-2}\log \sigma^{-1}$ and $\delta''\gg
\frac{\alpha^4\delta }{d^3\log\sigma^{-1}}$.

\end{proposition}

{\it Proof: }
We essentially follow the argument of Sanders \cite{S06}; however, some
care must be taken to get the right quantitative
version.
Replacing $f$ by $f(\cdot+b)$ if necessary, we may assume that $b=0$.
Let $c_0$ be a small enough constant which will be fixed
later.  By \cite{GT-inverse}, Lemma 8.2, we can find $\delta'$ such
that
\begin{equation}\label{zz.e5}
\delta'\in(c_0\alpha^2\delta d^{-1}, 2c_0\alpha^2\delta d^{-1})
\end{equation}
and that the set $B'$ defined in (i) is regular.  Suppose that
(\ref{pospropbohr}) fails for this choice of $\delta'$; we have to
prove that this implies (ii).

The failure of (\ref{pospropbohr}) means that we can find a set
$S\subset B'\cap \{x:(f*f)(x)<\frac{\alpha^2}{2}|B|\}$ such that
$|S|=\sigma|B'|$. Let $g=f-\alpha B$ be the ``balanced function" of
$f$.  We first claim that
\begin{equation}\label{convestimate}\frac{1}{|B||B'|}\sum_{x\in S}g\ast
g(x)\leq -\frac{\alpha^2\sigma}{2}+O(d\delta'\delta^{-1}\sigma).
\end{equation}
To prove this, we write
\begin{align*}
 \frac{1}{|B||B'|}\sum_{x\in S}(g\ast g)(x)
       =\frac{1}{|B||B'|}\big(\sum_{x\in S}(f\ast f)(x)-2\alpha\sum_{x\in S}(B\ast f)(x)
       +\alpha^2\sum_{x\in S}(B\ast B)(x)\big).
\end{align*}
The first term obeys
\begin{equation}\label{zz.e1}
\frac{1}{|B||B'|}\sum_{x\in S}(f\ast f)(x)
    \leq \frac{\alpha^2|B|}{2|B||B'|}|S|
    =\frac{\alpha^2\sigma}{2},\end{equation}
by the choice of $S$. The second term is estimated as in \cite{S06}.
By \cite{S06}, Corollary 3.4, we have for $x\in B'$
$$|f\ast \frac{B}{|B|}(x)-f\ast \frac{B}{|B|}(0)|\ll d\delta'\delta^{-1}.$$
But $f\ast \frac{B}{|B|}(0)= \alpha$, so that $f\ast
\frac{B}{|B|}(x) =\alpha+O(d\delta'\delta^{-1})$ for $x\in B'$.
Hence
\begin{equation}\label{zz.e2}
\frac{1}{|B'|}\sum_{x\in S}\frac{B}{|B|}\ast f(x)=
        \frac{|S|}{|B'|}(\alpha+ O(d\delta'\delta^{-1}))
    =\alpha\sigma+O(d\delta'\delta^{-1}\sigma),\end{equation}
Finally, we trivially have $B\ast B(x)\leq |B|$ for all $x$, hence
 \begin{equation}\label{zz.e3}
 \frac{1}{|B||B'|} \sum_{x\in S}B\ast B(x)
   \leq \sigma+O(d\delta'\delta^{-1}\sigma).\end{equation}
Combining (\ref{zz.e1}), (\ref{zz.e2}), (\ref{zz.e3}), we get
(\ref{convestimate}).

We now convert this to a Fourier analytic statement.  We have
\begin{align*}\sum_{x\in S}g\ast g(x) &=\sum_{x\in \zz}g\ast
        g(x)S(x)\\
    &=N\sum_{\xi\in\zz}\widehat{g\ast g}(\xi)\widehat{S}(\xi)\\
    &=N^2\sum_{\xi\in\zz}|\widehat{g}(\xi)|^{2}\widehat{S}(\xi).\end{align*}
Hence, by the triangle inequality, (\ref{convestimate}) implies that
\begin{equation}\label{fourierestimate}\frac{N^2}{|B||B'|}
\sum_{\xi}
        |\widehat{g}(\xi)|^2|\widehat{S}(\xi)|\geq
        \frac{\alpha^2\sigma}{2}+O(d\delta'\delta^{-1}\sigma).\end{equation}
Define $$\mathcal{L}:=\{\xi\in\zz:|\widehat{S}(\xi)|\geq
    \frac{\alpha\sigma|B'|}{4N}\}.$$  We claim that the main contribution to the
sum in (\ref{fourierestimate}) comes from $\mathcal{L}$.  In fact
\begin{align*}
\frac{N^2}{|B||B'|}\sum_{\xi\not\in\mathcal{L}}|\widehat{g}(\xi)|^2|\widehat{S}(\xi)|
    &\leq\frac{\alpha\sigma N}{4|B|}\sum_{\xi\not\in\mathcal{L}}|\widehat{g}(\xi)|^2\\
    &\leq\frac{\alpha\sigma N}{4|B|}\sum_{\xi\in\zz}|\widehat{g}(\xi)|^2\\
    &=\frac{\alpha\sigma}{4|B|}\sum_{x\in\zz}|g(x)|^2\\
    &=\frac{\alpha\sigma}{4|B|}\sum_{x\in\zz}|f(x)-\alpha B(x)|^2\\
    &=\frac{\alpha\sigma}{4|B|}\sum_{x\in\zz}f(x)^2- 2\frac{\alpha^2\sigma}{4|B|}\sum_{x\in\zz}f(x)B(x) +
        \frac{\alpha^3\sigma}{4|B|}\sum_{x\in\zz}B(x)^2\\
    &\leq \frac{\alpha^2\sigma}{4}-\frac{2\alpha^3\sigma}{4}+\frac{\alpha^3\sigma}{4}\\
    &=\frac{\alpha\sigma}{4}(\alpha-\alpha^2)\\
    &\leq \frac{\alpha^2\sigma}{4}\end{align*}

Hence
    $$\frac{N^2}{|B||B'|}\sum_{\xi\in\mathcal{L}}|\widehat{g}(\xi)|^2|\widehat{S}(\xi)|
        \geq \frac{\alpha^2\sigma}{4}+O(d\delta'\delta^{-1}).$$
Since $\frac{N}{|B'|}|\widehat{S}(\xi)|$ is trivially bounded by
$\sigma$, we have
\begin{equation}\label{zz.e4}
\frac{N}{|B|}\sum_{\xi\in\mathcal{L}}|\widehat{g}(\xi)|^2\geq
        \frac{\alpha^2}{4}+O(d\delta'\delta^{-1}).
\end{equation}

We now apply the localized version of Chang's theorem proved in
\cite{S06} (Proposition 4.2) to $S\subset B'$, with
$\epsilon=\alpha/4$ and $\eta=1/2$.  We conclude that there is a set
$\Lambda\subset\zz$ and a $\delta_{0}''>0$ such that
    $$|\Lambda|\ll\frac{2^4}{\alpha^2}\log\sigma^{-1},$$
    $$\delta_{0}''\gg \frac{\delta'\alpha^{2}4}{d^2\log\sigma^{-1}},$$
and
$$
\mathcal{L}\subset\{\xi\in\zz:|1-e^{-2\pi ix\xi/N}|\leq 1/2\
\forall\
        x\in B(\Gamma\cup\Lambda, \delta_{0}'')\}.
$$
Choose $\delta''\in(\delta_{0}'',2\delta''_0)$ such that
$B'':=B(\Gamma\cup\Lambda, \delta'')$ is regular.  Note that this
together with (\ref{zz.e5}) implies that $\delta''$ obeys the
condition in (ii).  We may also assume that $\delta''<\delta'$.  Our
goal is to get the $L^2$ density increment as in
(\ref{densityincrement}) on a translate of $B''$.

By the definition of $\mathcal{L}$, we have
$\frac{N}{|B''|}|\widehat{B''}(\xi)|\geq 1/2$ for all
$\xi\in\mathcal{L}$.  Hence
    $$\frac{N^3}{|B||B''|^2}\sum_{\xi\in\mathcal{L}}|\widehat{g}(\xi)
        |^2|\widehat{B''}(\xi)|^2\geq\frac{\alpha^2}{16}+O(d\delta'\delta^{-1}).$$
Again using Plancherel's identity and the convolution identity we
have
\begin{align*}\alpha^2\big(\frac{1}{16}+O(\alpha^{-2}d\delta'\delta^{-1})\big)&
        \leq\frac{N^3}{|B||B''|^2}\sum_{\xi\in\zz}|\widehat{g}(\xi)
        |^2|\widehat{B''}(\xi)|^2\\
    &=\frac{N^3}{|B||B''|^2}\sum_{\xi\in\zz}|N^{-1}\widehat{g\ast B''}(\xi)|^2\\
    &=\frac{1}{|B||B''|^2}\sum_{x\in\zz}|g\ast
        B''(x)|^2.\end{align*}
We now apply Lemma 5.2 from \cite{S06} and conclude that
\begin{align*}
\frac{1}{|B''|}\sup_{x\in\zz}|f\ast B''(x)|&\geq
        \alpha\big(1+2^{-4}+O(\alpha^{-2}d\delta'\delta^{-1})\big)+O(d\delta''\delta^{-1})\\
    &\geq\alpha\big(1+2^{-4}\big)+O(d\alpha^{-1}\delta'\delta^{-1}).\end{align*}
We now let the constant $c_0$ in (\ref{zz.e5}) be small
enough, so that the error term is bounded by $ \alpha 2^{-5}$.  The
conclusion (ii) follows if we choose $b''$ to maximize $|f\ast B''(b'')|$.

\begin{proposition}\label{f1estimate} Let $f:\zz\rightarrow [0,1]$ be defined
such that $$\ee_{x\in\zz}f(x)=\alpha>0.$$  Let $\sigma\in(0,1]$.
Then there exist Bohr sets $B:=b+B(\Gamma,\delta)$ and
$B':=b+B(\Gamma,\delta')$ such that
$$\Big|\{x\in B':\ f*f(x)\geq
\frac{\alpha^2}{2}|B|\}\Big|>(1-\sigma)|B'|,$$ and
$$\delta'\gg\frac{\alpha^2\delta}{|\Gamma|},$$
$$\delta\gg\big(\frac{\alpha}{\log(\sigma^{-1})}\big)^{C\log(\alpha^{-1})},$$
and $$|\Gamma|\ll\alpha^{-2}\log(\sigma^{-1}).$$

\end{proposition}

{\it Proof of Proposition \ref{f1estimate}:} We construct the Bohr
sets $B$ and $B'$ by iterating Proposition \ref{iteration}.
Let $\Gamma_0:=\{0\}$, and pick $\delta_0\gg 1$ so
that $B(\Gamma_0,\delta_0)$ is regular.  Define
$\alpha_0:=\alpha$. Averaging over translates of $B(\Gamma_0,\delta_0)$,
we see that there is a $b_0$ such that $\ee(f|B_0)\geq\alpha_0$ for
$B_0= b_0+B(\Gamma_0,\delta_0)$.
By Proposition \ref{iteration}, one of the following must hold:

(i) There is a $\delta_0'\gg
\frac{\alpha_0^2\delta_0}{|\Gamma_0|}$ such that
$B_0':=b_0+B(\Gamma_0,\delta_0')$ is regular and
\begin{equation}\label{pospropbohr2} \Big|\{x\in B_0':(f*f)(x)\geq
\frac{\alpha_0^2}{2}|B_0|\} \Big|\geq
(1-\sigma)|B_0'|,\end{equation}

(ii) There is a regular Bohr set
$B_1:=b_1+B(\Gamma_0\cup\Lambda_0,\delta_1)$ such that
\begin{equation}\label{densityincrement2}
\ee(f|B_1)\geq\alpha_0(1+2^{-5}) ,\end{equation}
where $|\Gamma_0|\ll\alpha_0^{-2}\log(\sigma^{-1})$ and $\delta_1\gg
\frac{\alpha_0^4\delta_0}{|\Gamma_0|^3\log(\sigma^{-1})}$.

If (i) holds, we let $B'=B'_0$ and we are done.  If on the other hand (ii)
holds, we repeat the procedure with $B_0$ replaced by $B_1$, and continue
by induction.
If we have not satisfied (i) by the end of the
$k$th step, we have found a regular Bohr set
$B_k:=b_k+B(\Gamma_k,\delta_k)$ such that
$$\ee(f|B_k)=\alpha_k|B_k|,$$
where
    \begin{equation}\label{alpha} \alpha_k\geq\alpha_{k-1}(1+2^{-5}),
        \end{equation}
    \begin{equation}\label{delta}
    \delta_k\gg\frac{\alpha_{k-1}^4\delta_{k-1}}{|\Gamma_{k-1}|^3\log(\sigma^{-1})},
        \end{equation}
and
    \begin{equation}\label{gamma-diff}
     |\Gamma_k|-|\Gamma_{k-1}| \ll\alpha_{k-1}\log(\sigma^{-1}).
        \end{equation}

The iteration must terminate (upon reaching density 1 on a
large enough Bohr set) after at most $$k\ll\log(\alpha^{-1})$$
steps, since from (\ref{alpha}) we have
$$\alpha_k^2\geq\alpha^2(1+2^{-5})^{k-1}.$$
By (\ref{gamma-diff}) we have
\begin{align*}
    |\Gamma_k|&\ll\alpha_{k-1}^{-2}\log(\sigma{-1})+
    \alpha_{k-2}^{-2}\log(\sigma{-1}) +\cdots
    +\alpha_{0}^{-2}\log(\sigma{-1})\\
    &\leq
    \alpha^{-2}\log(\sigma^{-1})\sum_{j=0}^{\infty}(1+2^{-5})^{-j}\ll
    \alpha^{-2}\log(\sigma^{-1}).\end{align*}
Finally, using our bounds for $\alpha_k$ and $|\Gamma_k|$, we have
$$\delta_k\gg\big(\frac{\alpha}{\log(\sigma^{-1})}\big)^{C\log(\alpha^{-1})},$$
for some absolute constant $C>0$.


\section{The restriction argument}\label{sec-random}

Assume that the hypotheses of Theorem \ref{sumsets-theorem} hold.
We need to show that if $f_1,f_2$ are as Lemma
\ref{lemma-decomposition} and $B,B'$ are the Bohr sets chosen in
Proposition \ref{f1estimate}, then (\ref{e-random}) holds, i.e.
\begin{equation}\label{e-random-bis}
\Big|\{x\in B':\ |f_2*f_i(x)|\geq \frac{\alpha^2}{10}|B|\}\Big|\leq
\sigma|B'|,\ i=1,2.
\end{equation}
It suffices to prove that
\begin{equation}\label{ww.e1}
\|f_i*f_2\|^2_{L^2(B')}\leq\frac{\alpha^4}{200}\sigma |B|^2|B'|.
\end{equation}
We have
\begin{align*}
\|f_i*f_2\|^2_{L^2(B')}&=\sum_{x\in B'}(f_i*f_2)^2(x) =\sum_{x\in
B'}\Big( \sum_{y}f_i(y)f_2(x-y)\Big)\Big(
\sum_{z}f_i(z)f_2(x-z)\Big)\\
&=\sum_{x,y,z,u,v} B'(x)f_i(y)f_2(z)\frac{1}{N}\sum_{\xi} e^{-2\pi
i(y+z-x)\xi/N}
\\
&\qquad\cdot f_i(u)f_2(v)\frac{1}{N}\sum_{\eta} e^{-2\pi
i(u+v-x)\eta/N}
\\
&=N^3\sum_{\xi,\eta}\widehat{B'}(-\eta-\xi)\widehat{f_i}(\xi)\widehat{f_2}(\xi)
\widehat{f_i}(\eta)\widehat{f_2}(\eta)
\\
&=N^3\sum_{\xi}(\widehat{B'}*\widehat{f_i}\widehat{f_2})(-\xi)
\widehat{f_i}(\xi)\widehat{f_2}(\xi).
\end{align*}
By H\"older's inequality,
\begin{equation}\label{vv.e1}
\|f_i*f_2\|^2_{L^2(B')} \leq
N^3\|\widehat{B'}*\widehat{f_i}\widehat{f_2}\|_{10}\,
\|\widehat{f_i}\widehat{f_2}\|_{10/9}.
\end{equation}
Applying Young's inequality, we get
\begin{equation}\label{vv.e2}
\|\widehat{B'}*\widehat{f_i}\widehat{f_2}\|_{10} \leq
\|\widehat{B'}\|_5\,\|\widehat{f_i}\widehat{f_2}\|_{10/9}.
\end{equation}
Furthermore,
\begin{align*}
\|\widehat{f_i}\widehat{f_2}\|_{10/9}^{10/9}
&\leq
\|\widehat{f_2}\|_\infty^{1/9}\,\sum_\xi|\widehat{f_2}(\xi)|\,|\widehat{f_i}(\xi)|^{10/9}
\\
&\leq \|\widehat{f_2}\|_\infty^{1/9}\,\|\widehat{f_2}\|_{19/9}\,
\|\widehat{f_i}(\xi)\|_{19/9}^{10/9},
\end{align*}
where at the last step we used H\"older's inequality again.
Plugging this together with (\ref{vv.e2}) in (\ref{vv.e1}), we see
that
\begin{align*}
\|f_i*f_2\|^2_{L^2(B')} & \leq
N^3\|\widehat{B'}\|_5\,\|\widehat{f_i}\widehat{f_2}\|_{10/9}^2
\\
& \leq N^3\|\widehat{B'}\|_5\,\Big(
\|\widehat{f_2}\|_\infty^{1/9}\,\|\widehat{f_2}\|_{19/9}\,
\|\widehat{f_i}\|_{19/9}^{10/9}\Big)^{9/5}
\\
&\leq N^3\|\widehat{B'}\|_5\,
\|\widehat{f_2}\|_\infty^{1/5}\,\|\widehat{f_2}\|_{19/9}^{9/5}\,
\|\widehat{f_i}\|_{19/9}^{2}.
\\
\end{align*}

By Plancherel's theorem and Lemma \ref{lemma-decomposition}(iv), we
have
$$\|\widehat{f_i}\|_2^{2}\leq\|\widehat{f}\|_2^2=
N^{-1}\|f\|^2\ll \alpha p^{-1}=\alpha N^\theta.$$ Since
$\theta<1/20$, it follows from Lemma \ref{lemma-restriction} that
$$\|\widehat{f}\|_{19/9}=O(1)\hbox{ and }\|\widehat{f_i}\|_{19/9}=O(1),
\ i=1,2.$$ By Lemma \ref{lemma-decomposition}(iii), we have
$$\|\widehat{f_2}\|_\infty\leq C\epsilon_0.$$
Finally,
$$
\|\widehat{B'}\|_5^5\leq
\|\widehat{B'}\|_\infty^3\,\|\widehat{B'}\|_2^2 \leq
\frac{|B'|^3}{N^3}\|\widehat{B'}\|_2^2 = \frac{|B'|^4}{N^4}.
$$
Combining these estimates, we get
\begin{equation}\label{vv.e6}
\|f_i*f_2\|^2_{L^2(B')} \ll N^3\epsilon_0^{1/9}
\frac{|B'|^{4/5}}{N^{4/5}}.
\end{equation}
We need the right side of this to be smaller than
$\frac{\alpha^4}{200}\sigma |B|^2|B'|$, i.e. we need to have
\begin{equation}\label{vv.e7}
\epsilon_0^{1/9}\leq c\alpha^4\sigma
\frac{|B|^2}{N^2}\frac{|B'|^{1/5}}{N^{1/5}} =c\alpha^4\sigma
\pp(B)^2\pp(B')^{1/5}.
\end{equation}
But by Lemma \ref{lemma-bohrsize} and (\ref{pp.e4})--(\ref{pp.e6}),
$\pp(B)$ and $\pp(B')$ are bounded from below by
$$\pp(B)\geq\pp(B')\gg(c\delta'')^{|\Gamma|}
\gg\Big(\frac{c\alpha}{\log k}\Big)^{c\alpha^{-2}\log\frac{1}{\alpha}\log k},
$$
where we plugged in $\sigma=(16k)^{-1}$.
Hence (\ref{vv.e7}) holds if
\begin{equation}\label{pp.e7}
\epsilon_0\ll\alpha^{28}k^{-9}
\Big(\frac{c\alpha}{\log k}\Big)^{c\alpha^{-2}\log\frac{1}{\alpha}\log k}.
\end{equation}
A short calculation shows that (\ref{pp.e8}) is sufficient to guarantee that
(\ref{pp.e7}) is satisfied.


\section{Proof of Proposition \ref{sumsets-size}}\label{sec-last}
\init

Let $0<\sigma<(\alpha-\beta)/10$.
Define $\nu,f,f_1,f_2$ as in Section \ref{long-section}, except that instead of
(\ref{pp.e3}) we will require
\begin{equation}\label{ff.e1}
\|f_1\|_\infty\leq 1+\sigma,
\end{equation}
which holds for large enough $N$ (depending on $\sigma$ and on the $\epsilon_0$
in the definition of $f_i$) by the same argument as in Section
\ref{sarkozy}.

It clearly suffices to prove that
\begin{equation}\label{ff.e2}
\Big|\{x\in\zz:\ f*f(x)>0\}\Big|\geq (\alpha-10\sigma)N.
\end{equation}
Indeed, (\ref{ff.e2}) shows that the sumset $A+A$ in $\zz$ has size at least $\beta N$, hence
so does the sumset $A+A$ in $\mathbb{Z}$.

We first claim that if $N$ is large enough, then
\begin{equation}\label{ff.e3}
\Big|\{x\in \zz:\ f_1*f_1(x)\geq \sigma\alpha N\}\Big|\geq(\alpha-3\sigma)N.
\end{equation}
To see this, we first note that
\begin{equation}\label{ff.e4}
\|f_1*f_1\|_1=\|f_1\|_1^2=\alpha^2N^2(1+O(N^{-1/5})).
\end{equation}
On the other hand, if (\ref{ff.e3}) failed, we would have
$$
\|f_1*f_1\|_1\leq \sigma\alpha N\cdot N+\alpha N(1+\sigma+ O(N^{-1/5}))
\cdot(\alpha-3\sigma)N$$
$$
=\alpha^2 N^2(1+ O(N^{-1/5}))-\sigma\alpha N^2,$$
which contradicts (\ref{ff.e4}).  This proves (\ref{ff.e3}).

The proof of (\ref{ff.e2}) will be complete if we can show that
\begin{equation}\label{ff.e5}
\Big|\{x\in \zz:|\ f_i*f_2(x)|\geq \frac{\sigma\alpha}{10} N\}\Big|\leq\sigma N.
\end{equation}
To this end, we repeat the argument in Section \ref{sec-random}.  It suffices
to prove that
\begin{equation}\label{ff.e6}
\|f_i*f_2\|^2_{2}\leq\frac{\sigma^2\alpha^2}{200}\sigma N^3.
\end{equation}
As in Section \ref{sec-random} (with $B=B'=\zz$), we have
\begin{equation}\label{ff.e7}
\|f_i*f_2\|^2_2 \ll \epsilon_0^{1/9} N^3,
\end{equation}
and the right side is smaller than the right side of (\ref{ff.e6}) if
$\epsilon_0\ll \sigma^{27}\alpha^{18}$, with a small enough implicit constant.
Thus (\ref{ff.e5}) holds for large enough $N$ if $\epsilon_0$ was chosen
small enough.


\section{Acknowledgements}

The authors were supported in part by an NSERC Discovery Grant.  We are grateful
to Ben Green for suggesting the feasibility of Theorem \ref{sarkozyextension},
and to Ernie Croot and Mihalis Kolountzakis for helpful discussions
and suggestions.

\noindent{\sc Department of Mathematics, University of British Columbia, Vancouver,
B.C. V6T 1Z2, Canada}

\noindent{\it melamel@math.ubc.ca, ilaba@math.ubc.ca}


\end{document}